\documentclass[a4paper,14pt]{article}
\usepackage[14pt]{extsizes}
\usepackage{mathtext}                    
\usepackage{cmap}                        
\usepackage[cp1251]{inputenc}            
\usepackage[english, russian]{babel}     
\usepackage[left=2.5cm,right=1cm,top=1cm,bottom=2cm]{geometry} 
\usepackage{amssymb}
\usepackage{amsmath}

\begin{document}

\newtheorem{theorem}{Theorem}[section]
\newtheorem{corollary}[theorem]{Corollary}
\newtheorem{lemma}[theorem]{Lemma}
\newtheorem{proposition}[theorem]{Proposition}

\newtheorem{definition}[theorem]{Definition}
\newtheorem{example}[theorem]{Example}

\newtheorem{remark}[theorem]{Remark}


 \centerline{\textbf{\textsc{The Markov-Stieltjes  transform as an operator}}}
\vspace{1cm}

 \centerline{A.R. Mirotin, I.S. Kovalyova}

\vspace{1cm}


We prove that the Markov-Stieltjes  transform  is a bounded non compact Hankel operator on
Hardy space $H^p$ with  Hilbert matrix with respect to  the standard Schauder  basis of $H^p$ and
a bounded non compact operator on Lebesgue space $L^p[0,1]$ for $p\in(1,\infty)$ and obtain
estimates for its norm in this spaces. It is shown  that the Markov-Stieltjes  transform on
$L^2(0,1)$ is unitary equivalent to the Markov-Stieltjes  transform on  $H^2$.
Inverse formulas and operational properties for this transform are obtained.

\vspace{1cm}

\section{Introduction}
 Much work has been done during  last years on the theory of  integral transforms of functions of one real variable and in particular on convolution and inversion theorems for such transforms and their applications to integral equations (see, e.g., \cite{Ya1} -- \cite{KM} and the bibliography cited therein). This paper is devoted to the Markov-Stieltjes  transform $S$ of functions on $(0,1)$.
The last  transform   was introduced in \cite[Chapter 6]{Mir} as a special case of the Stieltjes  transform of measures on general semigroups. The terminology goes back to approximation theory, see, e.g., \cite[p. 14]{And}, \cite{And2}, \cite{MGU}.   We give inverse formulas for this transform and formulate its operational properties. The main goal of this paper is to study the Markov-Stieltjes  transform as an operator on Hardy spaces $H^p$ for $p\in(1,\infty]$ and  Lebesgue spaces $L^p(0,1)$ for $p\in(1,\infty)$. We prove that $S$ is
 a bounded non compact Hankel operator on  Hardy space $H^p$ with  Hilbert matrix with respect to  the standard Schauder  basis of $H^p$ for $p\in(1,\infty)$, a bounded non compact operator from  $H^\infty$ to $BMOA$, and  a bounded  operator in  $\ell_A^p$ and give estimates for the norm of $S$ in this cases. We show also that $S$ is
 a bounded non compact operator  on Lebesgue space $L^p(0,1)$ for $p\in(1,\infty)$ and obtain  estimates for its norm in this spaces, too.
 It is shown also that the Markov-Stieltjes  transform on  $L^2(0,1)$ is unitary equivalent to the Markov-Stieltjes  transform on  $H^2$. As a corollary the norm and the spectrum of $S$ as an operator on
$L^2(0,1)$ are obtained.

\begin{definition}\cite[Chapter 6]{Mir},  \cite{KM}.
The Markov-Stieltjes  transform of a  function $f\in L^1(0,1)$  is  defined by the formula
\begin{equation}
Sf(z):=\int_0^1 \frac{f(t)}{1-tz}dt
\end{equation}
\end{definition}
\noindent(we write also $S_{t\to z}\{f(t)\}$ instead of $Sf(z)$). Obviously, for $z\notin [1,\infty)$ this Lebesgue integral exists  and represents an analytic function $f^*$ in the domain $\mathbb{C}\setminus[1,\infty)$. For $z\in [1,\infty)$ the integral in (1) is understood as a Cauchy principle value integral, i.e.,
\begin{equation}
Sf(z):=V.P.\int_0^1 \frac{f(t)}{1-tz}dt:=\lim_{\varepsilon\to 0+}\int\limits_{\{t\in (0,1):|t-1/z|>\varepsilon\}} \frac{f(t)}{1-tz}dt.
\end{equation}
The limit in the right hand side of (2)  exists for almost all $z\in [1,\infty)$. In fact,
$Sf(z)=(\pi/z)(Hf_1)(1/z)$, where $Hf_1$ stands for the Hilbert transform of the function
$f_1(x):=f(x)$ for $x\in (0,1)$ and $f_1(x):=0$ otherwise. So, the application of Loomis Theorem (see, e.g., \cite[p. 239]{King})  proves the assertion.

The following example shows that in general $\mathbb{C}\setminus[1,\infty)$ is the domain of holomorphy of $Sf$.

\begin{example}
Putting $t=x^2/(1+x^2)$ it is easy to verify that $S_{t\to z}\{(t(1-t))^{-1/2}\}$ equals to $\pi(1-z)^{-1/2}$ for $z\notin [1,\infty)$ and equals to zero otherwise.
\end{example}

As was mentioned above the study of  Markov-Stieltjes  transform as a
function is important in approximation theory (\cite[p. 226]{And}, \cite{MGU}).


This transform is useful in solving some singular integral equations, too. In \cite{KM} for $f\in L^p(0,1), g\in L^q(0,1)\  (1<p, q<\infty, 1/p+1/q<1)$ the following  binary operation was considered (this operation was introduced for the first time in   \cite[p. 220, formula (24.38)]{Ya1})

\[
f\circledast g(t)=tf(t)\int\limits_0^1\frac{f(u)}{t-u}du+tg(t)\int\limits_0^1\frac{g(u)}{t-u}du,
\]
where the integrals are understood as their Cauchy   principal values, and using methods developed by H.M. Srivastava and Vu Kim Tuan in \cite{ST}  a  convolution theorem for Markov-Stieltjes  transform in the form
\[
S(f\circledast g)=(Sf)\cdot (Sg)
\]
was proved. Arguing as in \cite{ST} it was also shown in \cite{KM} that the equation
\[
x(t)+\lambda\int\limits_0^1\frac{x(u)}{t-u}du=g(t)\ (\lambda\ne 0)
\]
where $g$ is prescribed and $x$ is an unknown function to be determined has (for appropriate $g$) the unique solution
\[
x(u)=\cos(\alpha\pi)S_{s\to u}^{-1}\left\{(1-s)^\alpha S_{t\to s}\left\{ \frac{g(t)t^\alpha}{(1-t)^\alpha}\right\}\right\},
\]
$\alpha$  being a (unique) root
of the  equation $\tan(\alpha\pi)=\lambda\pi,\ 0<{\rm Re}\alpha<1$ (for inversion formulas for the Markov-Stieltjes  transform see the following section).

\section{Inversion formulas}

A complex inversion formula for the Markov-Stieltjes  transform looks as follows.

\begin{theorem}
Let $f\in L^1(0,1), 0<t<1$, and $f(t\pm 0)$ exist. If $f^*=Sf$, then
\[
\frac{f(t+0)+f(t-0)}{2}=\frac{1}{2\pi i}\lim\limits_{\eta\to 0+}\left(\frac{1}{t-i\eta}f^*\left(\frac{1}{t-i\eta}\right)-\frac{1}{t+i\eta}f^*\left(\frac{1}{t+i\eta}\right)\right).
\]
\end{theorem}

Proof.
It is easy  to verify that
\[
\frac{1}{2\pi i}\left(\frac{1}{t-i\eta}f^*\left(\frac{1}{t-i\eta}\right)-\frac{1}{t+i\eta}f^*\left(\frac{1}{t+i\eta}\right)\right)=
\frac{1}{\pi}\int\limits_0^1\frac{\eta}{(t-s)^2+\eta^2}f(s)ds.
\]
Application of \cite[p. 338, Lemma 7.2]{Widder} completes the proof.

Now we shall formulate also a real inversion formula for the Markov-Stieltjes  transform.
\begin{theorem} \cite{KM}.
Let $f\in L^p(0,1), 1<p<\infty$. The Markov-Stieltjes  transform $f^*(x)=Sf(x)$ exists for a.e. $x\in\mathbb{ R}$ and
\[
f(t)=\frac{1}{\pi^2}V.P.\int\limits_{-\infty}^\infty\frac{f^*(x)}{1-tx}dx.
\]
\end{theorem}

Proof.
It follows, e.g., from the above mentioned  equality $Sf(z)=(\pi/z)(Hf_1)(1/z)$ and the inversion formula for the Hilbert transform (see \cite{KM} for details).

\section{Operational properties of Markov-Stieltjes  transforms}

The following properties hold for Markov-Stieltjes  transform  (cf., e. g., \cite[p. 394]{DB}).

If  $f^*=Sf$, then

1) $$S_{t\to z}\{f(1-t)\}=\frac{1}{1-z}f^*\left(\frac{z}{z-1}\right);$$

2) $$S_{t\to z}\{f(at)\}=\frac{1}{a}f^*\left(\frac{z}{a}\right)\ (a>1);$$

3) $$S_{t\to z}\{tf(t)\}=\frac{1}{z}\left(f^*(z)-\int_0^1f(t)dt\right)\ (f \in L^1(0,1)).$$

 In particular, $S_{t\to z}\{tf(t)\}=\frac{1}{z}f^*(z) \mbox{ if }  \int_0^1f(t)dt=0;$

4) $$S_{t\to z}\left\{\frac{f(t)}{t+a}\right\}=\frac{z}{1+az}f^*(z)+\frac{1}{1+az}\int_0^1\frac{f(t)}{t+a}dt\ \left(\frac{f(t)}{t+a}\in L^1(0,1)\right).$$

In particular, $S_{t\to z}\left\{\frac{f(t)}{t+a}\right\}=\frac{z}{1+az}f^*(z) \mbox{ if } \int_0^1\frac{f(t)}{t+a}dt=0;$

5) $$S_{t\to z}\left\{\frac{d}{dt}f(t)\right\}=-\frac{d}{dz}f^*(z)+\frac {f(1)}{1-z}-f(0)  \mbox{ if } f\in C^1[0,1], \ z\notin [1,\infty).$$

In particular, $S_{t\to z}\{\frac{d}{dt}f(t)\}=-\frac{d}{dz}f^*(z) \mbox{ if in addition } f(0)=f(1)=0$;

6) $$S_{t\to z}\left\{\int_0^tf(t)dt\right\}=-\int_0^zf^*(z)dz-\left(\int_0^1f(t)dt\right)\log(1-z)+\int_0^1(1-t)f(t)dt.$$

In particular, $S_{t\to z}\{\int_0^tf(t)dt\}=-\int_0^zf^*(z)dz  \mbox{ if } \int_0^1f(t)dt=\int_0^1tf(t)dt=0$.

We omit simple proofs of this properties.

\section{Markov-Stieltjes  transform as an operator on Hardy spaces}

In this section we  identify the Hardy spaces $H^p(\mathbb{D})$  and $H^p(\mathbb{T})$ ($\mathbb{D}$ stands for the open unit disk and $\mathbb{T}$ for the unit sircle; see, e.g., \cite{Dur}) and frequently  use the notation $H^p$ for this space,  we denote also by
$\chi_n(z):=z^n\ (n\in \mathbb{Z}_+)$ the standard  (Schauder) basis of $H^p(\mathbb{D})$.

\begin{definition}
Following  \cite[p. 52]{Bott} for $b\in{L^\infty}(\mathbb{T})$ we define the Hankel operator $H(b)$ on $H^p(\mathbb{T}) (1<p<\infty)$ by
 $$
 H(b)\colon H^p \rightarrow H^p \colon f\mapsto PM(b)(I-P)Jf,
 $$
where
$$
P : \sum_{n=-N}^N f_n\chi_n \mapsto \sum_{n=0}^N f_n\chi_n,
$$
$$
M(b) : L^p(\mathbb{T})\rightarrow L^p(\mathbb{T}) : f\mapsto bf,
$$
$$
J : f(t)\mapsto \frac{1}{t}f\Bigl(\frac{1}{t}\Bigl)=\sum_{n\in \mathbb{Z}}f_n \chi_{-n-1}(t) (t\in\mathbb{T}),
$$
where $f=\sum_{n\in \mathbb{Z}}f_n \chi_n$.
\end{definition}
The function $b$ is called the symbol of the Hankel operator $H(b)$.

The next theorem describes the properties of $S$ as an operator on  $H^p(\mathbb{D})$ (this means that  $S$ is defined by the formula (1), where $f\in H^p(\mathbb{D}), z\in \mathbb{D}$).

\begin{theorem}
$\mathrm{1)}$ The Markov-Stieltjes  transform $S$ is a bounded non compact Hankel operator on  $H^p(\mathbb{D})\ (1<p<\infty)$ and has Hilbert matrix with respect to  the standard basis. Moreover, the following estimates hold
\begin{equation}
\pi\leq\| S\|_{H^p\rightarrow H^p}\leq \frac{\pi}{\sin \frac{\pi}{\max\{p,q\}}}.
\end{equation}
In particular, if $p=2$ then the norm and the essential norm of $S$ equal to $\pi$, and the spectrum and the essential spectrum
of $S$ equal to $[0,\pi]$.

$\mathrm{2)}$ The Markov-Stieltjes  transform $S$  is a bounded non compact Hankel operator from
$H^\infty(\mathbb{T})$ to $BMOA(\mathbb{T})$, and
$$
\|S\|_{H^\infty\rightarrow BMOA} \leq \pi\|P\|_{L^\infty\rightarrow BMOA}.
$$
\end{theorem}

 Proof.
  1) First note, that for  $f\in H^p(\mathbb{D})$ the  Fejer-Riesz inequality (see, e.g., \cite[Theorem 3.13]{Dur}) implies that
  \begin{equation}
\|f|(0,1)\|_{L^p(0,1)}\leq \pi^{1/p}\|f\|_{H^p}
\end{equation}
It follows that the restriction $f|(0,1)$ belongs to $L^p(0,1)$ and therefore the Lebesgue integral in (1) exists for all $p\in(1, \infty)$ and $z\in \mathbb{D}$.

Next, since for all $z\in \mathbb{D}$
\begin{equation}
(S\chi_n)(z)=\int_0^1 \frac{t^n}{1-tz}dt=\sum_{m=0}^{\infty}z^m\int_0^1 t^{n+m}dt=
\sum_{m=0}^{\infty}\frac{\chi_m(z)}{n+m+1},
\end{equation}
operator  $S$ is Hankel and has Hilbert matrix $\Gamma=(1/n+m+1)_{m,n=0}^\infty$ with respect to  the standard basis  $(\chi_n)_{n\in\mathbb{Z}_+}$ of $ H^p(\mathbb{D})$. Indeed, by  formula (5) and  Parceval's formula,
\[
\left\|S\chi_n-\sum_{m=0}^{M}\frac{\chi_m}{n+m+1}\right\|_{H^2}^2=\left\|\sum_{m=M+1}^{\infty}\frac{\chi_m}{n+m+1}\right\|_{H^2}^2=
\]
\[\sum_{m=M+1}^{\infty}\frac{1}{(n+m+1)^2}\to 0\ (M\to\infty).
\]
This implies that $S\chi_n=\sum_{m=0}^{\infty}\chi_m/(n+m+1)$ in the sense of $H^2$ and therefore  $\langle S\chi_j,\chi_k\rangle=1/(j+k+1)$ for
all $j, k\geq 0$ where $\langle \cdot,\cdot\rangle$ denotes the inner product in $H^2$.

Now we use the following remark  to the Nehari Theorem. If the $n$-th Fourier coefficient of a function $a\in H^p(\mathbb{D})$ equals to $a_n$ for $n\in \mathbb{N}$, and the operator $A$ in $H^p(\mathbb{T})$ satisfies $\langle A\chi_j,\chi_k\rangle=a_{j+k+1}$ for
all $j, k\geq 0$ then  $A$ is bounded on $H^p(\mathbb{T})$ if (and only if) $a\in BMO$, the space of functions of bounded mean oscillation on  $\mathbb{T}$   \cite[p. 55]{Bott}. But
\[
-\log(1-z)=\sum_{n=1}^\infty\frac{1}{n}\chi_n(z),
\]
 and  the function $a(z):=-\log(1-z)$ belongs to   $BMO$ ($a$ is analytic in $\mathbb{D}$ and its imaginary part  belongs to   $L^\infty(\mathbb{T})$, thus, $a$ has the form $f+\tilde{g}$, where $f,g\in L^\infty(\mathbb{T})$ and $\tilde{g}$ is the harmonic conjugate of $g$). It follows that $a\in H^p(\mathbb{D})$, since $BMO\subset L^p(\mathbb{T})$.
 Now by  the previous remark to the Nehari Theorem,    $S$ is bounded on $H^p(\mathbb{D})$.

Moreover, if $b$ is the symbol of   $S$  the Nehari Theorem \cite[Theorem 2.11]{Bott} implies
\[
{\rm dist}_{L^\infty}(b, \overline {H^\infty})\leq\left\| {S}\right\|_{H^p\rightarrow H^p}\leq c_p {\rm dist}_{L^\infty}(b, \overline {H^\infty}),
\]
where
$$
{\rm dist}_{L^\infty}(b, \overline {H^\infty})=\inf\{\|b-f\|_{L^\infty}:f\in \overline {H^\infty}\},\ c_p=\frac{1}{\sin \frac{\pi}{\max\{p,q\}}}
$$
(see, e.g., \cite[p. 32]{Bott}).
But the symbol of the Hankel operator  $S$ on $H^p$ does not depend of $p$ (see the proof of the Nehari Theorem in \cite{Bott}). So for $p=2$ we have
 \[
 \|S\|_{H^2\rightarrow H^2}={\rm dist}_{L^\infty}(b, \overline {H^\infty}).
 \]
On the other hand, it is known  that the norm and the essential norm of the Hankel operator on the space $H^2(\mathbb{T})$ with  Hilbert matrix with respect to  the standard  basis  equal to  $\pi$ (and therefore  ${\rm dist}_{L^\infty}(b, \overline {H^\infty})=\pi$; see, e.g., \cite[p. 36]{Pell}), and its spectrum  and   essential spectrum
 equal to $[0,\pi]$  (see, e.g., \cite[p. 575, Theorem 1.7]{Pell}).

To prove 1) it remains to show that $S$ is non compact on $H^p\ (1<p<\infty)$. For this we recall that the symbol  of the Hankel operator on the space $H^2(\mathbb{T})$ with    Hilbert matrix with respect to the standard  basis is  $b(e^{it})=ie^{-it}(\pi-t),\ 0\leq t < 2\pi$ (see, e.g., \cite[p. 6]{Pell}). According to the Hartman Theorem \cite[p. 80]{Bott} if the operator $S=H(b)$ is compact then  $b\in C(\mathbb{T})+\overline{H^\infty}$. But this inclusion contradicts the  Lindel$\ddot{\rm o}$f theorem  on one-sided  limits of $H^\infty$-functions (see, e.g., \cite[Corollary 5.3.5]{Nik}). This completes the proof of the first part of the theorem.

2) As it was shown above $Sf=PM(b)(I-P)Jf$ for $f\in H^p$.
Since $J : H^\infty \rightarrow (H^2)^\bot$,
it follows that  $Sf=PM(b)Jf$ for $f\in H^\infty$.
Moreover,
$$
J : H^\infty \rightarrow L^\infty, \|J\|_{H^\infty \rightarrow L^\infty}=1,
$$
and
$$
M(b) : L^\infty \rightarrow L^\infty, \|M(b)\|_{L^\infty \rightarrow L^\infty}=\|b\|_{L^\infty}=\pi.
$$
It is also known   (see, e.g., \cite[Theorem 8.3.10]{Zh}) that $P$ is a bounded operator from $L^\infty$ to $BMOA$. This implies that $S$ is a bounded operator from $H^\infty$ to $BMOA$ and
$$
\|S\|_{H^\infty \rightarrow BMOA} \leq \|P\|_{L^\infty \rightarrow BMOA}\|M(b)\|_{L^\infty \rightarrow L^\infty}\|J\|_{H^\infty \rightarrow L^\infty}=\pi\|P\|_{L^\infty \rightarrow BMOA}.
$$

 Finally, since $b\notin C(\mathbb{T})+\overline{H^\infty}$, the operator $S:H^\infty (\mathbb{T})\to BMOA(\mathbb{T})$ is non compact by the main result of the paper  \cite{CJY}.

Results like the previous theorem may have applications to approximation theory. Let $\mathbf{r}_n$ be the set of rational functions of order at most $n$ whose poles lies outside of $\mathbb{D}$, $X$ a Banach space of functions defined on some set $E\subseteq\mathbb{D}$ and $XR_n(f)=\inf_{r\in\mathbf{r}_n}\|f-r\|_X$ the degree of approximation of a function $f$ from $X$
by rationals from $\mathbf{r}_n$. By the triangle inequality,
\begin{equation}
|XR_n(f)-XR_n(g)|\leq\|f-g\|_{X}\ (f,g\in X).
\end{equation}

\begin{corollary}
The map $f\mapsto H^pR_n(Sf)$ is continuous on  $H^p$.
\end{corollary}

 Proof.
Indeed, the map  $g\mapsto H^pR_ng$ is continuous on  $H^p$ by the inequality (6).

The following corollary is a  generalization (for a Hausdorff moment problem) of a result  due to K. Zhu \cite[p. 372, Proposition 9]{Li}.

\begin{corollary}
Let $1<p\leq 2$. For $f\in H^p$ we let $Tf$ be the sequence $(c_n)$ defined by
\[
c_n=\int_0^1f(t)t^ndt,\ n\in \mathbb{Z}_+.
\]
Then $T$ is a bounded linear operator from $H^p$ to $\ell^q\ (1/p+1/q=1)$ and $\|T\|_{H^p\to \ell^q}\leq\pi/\sin\frac{\pi}{q}$.
\end{corollary}

 Proof.
Since
$$
Sf(z)=\int_0^1f(t)\sum_{n=0}^\infty(tz)^ndt =\sum_{n=0}^\infty c_nz^n,
$$
Theorem 6.1 from \cite{Dur} and Theorem 4.2 imply  that
$$
\|Tf\|_{\ell^q}\leq\|Sf\|_{H^p}\leq \frac{\pi}{\sin\frac{\pi}{q}}\|f\|_{H^p}.
$$

For the following corollary recall that  the Bergman space  $L^2_a$ consists of such  functions $f(z)=\sum_{m=0}^\infty f_mz^m$ that are holomorphic in $\mathbb{D}$ and
$$
\|f\|_{L^2_a}:=\left(\sum_{m=0}^\infty \frac{|f_m|^2}{m+1} \right)^{1/2}<\infty,
$$
 the sequence $\xi_n:=\sqrt{n+1}\chi_n$ forms an  orthonormal basis for  $L^2_a$ (see, e.g., \cite{Zh}).

\begin{corollary}
The Markov-Stieltjes  transform $S$ is an unbounded  densely defined operator on the Bergman space $L^2_a$.
\end{corollary}

 Proof.
Indeed, by the formula (5)
\[
S\xi_n(z)=\sum_{m=0}^\infty\frac{ \sqrt{n+1}}{m+n+1}z^m.
\]
So, the matrix $(a_{jk})$ of $S$ with respect to the basis $(\xi_j)$ is
\[
a_{jk}=\frac{\sqrt{k+1}}{\sqrt{j+1}(k+j+1)}.
\]
Since $\sum_{k=0}^\infty|a_{jk}|^2=\infty$, the operator $S$ is  unbounded. It remains to note that $H^2$ is dense in $L^2_a$.

\section{The Markov-Stieltjes  transform as an operator on Lebesgue spaces}

  In the following theorem $S$ denotes the Markov-Stieltjes  transform  on $L^p(0,1)\ (1<p<\infty)$. In other words, $S$ is defined by the formula (1), where $f\in L^p(0,1), z\in (0,1)$.

\begin{theorem}
$\mathrm{1)}$ The Markov-Stieltjes  transform  is a bounded non compact operator on  $L^p(0,1)\ (1<p<\infty)$. Moreover, the following estimates hold
$$
\frac{\pi}{\sin \frac{\pi}{p}}\leq \|S\|_{L^p\rightarrow L^p}\leq \pi{\cot \frac{\pi}{2\max\{p,q\}}}.
$$

$\mathrm{2)}$ The Markov-Stieltjes  transform  on $L^2 (0,1)$ is unitarily equivalent to   the Markov-Stieltjes  transform on  $H^2 (\mathbb{D})$. In particular, the norm and the essential norm of $S$ equal to $\pi$, and the spectrum and the essential spectrum
of $S$ equal to $[0,\pi]$.
\end{theorem}

 Proof.
We begin with the case  $1<p\leq 2$. As was mentioned in the Introduction, $ Sf(z)=y\pi Hf_1(y)$
 where   $z=1/y\ (y>0)$, $H$ stands for the Hilbert transform of functions on $\mathbb{R}$, and the function $f_1(t):=f(t)$ for $t\in(0,1)$ and $f(t):=0$ for $t\in \mathbb{R}\setminus(0,1)$
belongs to $L^p(\mathbb{R})$. Now the M. Riesz inequality for the Hilbert transform  implies for $1<p\leq2$ that
$$
\|Sf\|_{L^p(0,1)}=\Bigl(\int_0^1\Bigl|\frac{\pi}{z}Hf_1\Bigl(\frac{1}{z}\Bigl)\Bigl|^p dz\Bigl)^{1/p}=\Bigl(\int_1^\infty\frac{1}{y^2}|y\pi Hf_1(y)\bigl|^p dy\Bigl)^{1/p}=$$
$$=\pi\Bigl(\int_1^\infty\frac{1}{y^{2-p}}|Hf_1(y)\bigl|^p dy\Bigl)^{1/p}\leq\pi\Bigl(      \int_1^\infty|Hf_1(y)|^p dy\Bigl)^{1/p}\leq\pi\|Hf_1\|_{L^p(\mathbb{R})}\leq\pi A_p\|f\|_{L^p(0,1)}.
$$
Since (see, e.g., \cite{King})
$$
A_p=
\left\{
\begin{aligned}
\tan \frac{\pi}{2p},\ 1<p\leq 2
\\
\cot \frac{\pi}{2p},\ p>2\ \ \ \
\end{aligned}
\right.,
$$
we have
$$
\|S\|_{L^p\rightarrow L^p}\leq\pi\cot \frac{\pi}{2\max\{p,q\}}.
$$

 In the case $p>2$ using  standard duality arguments, this inequality, and H$\ddot{\rm o}$lder inequality we get (below for $f,g\in L^p(0,1)$ we put $\langle f, g\rangle :=\int_0^1 f\overline{g}dt$,  $A_q:=A_p,\ 1/p+1/q=1$)
$$
\|Sf\|_{L^p}=\sup\{\langle Sf, g\rangle : g\in L^q,\ \|g\|_{L^q}\leq 1\}=
\sup\{\langle S\overline{g}, \overline{f}\rangle : g\in L^q,\ \|g\|_{L^q}\leq 1\}\leq$$
$$\leq\sup\{\|S\overline{g}\|_{L^q} \|\overline{f}\|_{L^p} : g\in L^q,\ \|g\|_{L^q}\leq 1\}\leq$$
$$\leq\sup\{A_q\pi\|g\|_{L^q} \|f\|_{L^p} : g\in L^q,\ \|g\|_{L^q}\leq 1\}\leq A_q\pi\|f\|_{L^p}.
$$
This proves  the right-hand side of the desired  inequality.

To prove the left-hand side of this  inequality, consider the function
$$
f_\gamma (t):=\Biggl(\frac{t}{1-t}\Biggl)^\gamma, \gamma\in \Bigl(-\frac{1}{p}, \frac{1}{p}\Bigl).
$$
Then
$$
\|f_\gamma\|_{L^p}^p=\int_0^1\Biggl|\frac{t}{1-t}\Biggl|^{p\gamma}dt
=B(1+p\gamma, 1-p\gamma)=\frac{\pi p \gamma}{\sin\pi p\gamma}.
$$
Using \cite[Section 2.2.6, formula 5]{BMP} we have
$$
Sf_\gamma (t)=-\frac{\pi}{\sin\pi\gamma}\frac{1}{z}\Biggl(1-\frac{1}{(1-z)^\gamma}\Biggl).
$$
Therefore
$$
\|Sf_\gamma\|_{L^p}^p=\Bigl|\frac{\pi}{\sin\pi\gamma}\Bigl|^p\int_0^1\frac{1}{z^p}\Bigl|1-\frac{1}{(1-z)^\gamma}\Bigl|^pdz=
\Bigl(\frac{\pi}{\sin\pi\gamma}\Bigl)^p\int_0^1\Bigl(\frac{1-x^\gamma}{1-x}\Bigl)^p\frac{dx}{x^{p\gamma}}.
$$
Fix $0<\gamma_0<\frac{1}{p}$. For every
 $\varepsilon >0$ there exists such $\delta>0$ that $(1-x^{\gamma_0})^p(1-x)^{-p}>1-\varepsilon$  for all  $x\in(0, \delta)$. Then for $\gamma>\gamma_0$ and $x\in(0, \delta)$ we have
$$
\Bigl(\frac{1-x^{\gamma}}{1-x}\Bigl)^p>\Bigl(\frac{1-x^{\gamma_0}}{1-x}\Bigl)^p>1-\varepsilon.
$$

It follows that
$$
\Biggl(\frac{\|Sf_\gamma\|_{L^p}}{\|f_\gamma\|_{L^p}}\Biggl)^p=\frac{\sin\pi p\gamma}{\pi p \gamma}\Bigl(\frac{\pi}{\sin\pi\gamma}\Bigl)^p\int_0^1\Bigl(\frac{1-x^\gamma}{1-x}\Bigl)^p\frac{dx}{x^{p\gamma}}\geq
$$
$$
\geq\frac{\sin\pi p\gamma}{\pi p \gamma}\Bigl(\frac{\pi}{\sin\pi\gamma}\Bigl)^p\int_0^\delta (1-\varepsilon)\frac{dx}{x^{p\gamma}}=\frac{\sin\pi p\gamma}{\pi p \gamma}\Bigl(\frac{\pi}{\sin\pi\gamma}\Bigl)^p\frac{\delta^{1-p \gamma}}{1-p\gamma} (1-\varepsilon).
$$
Since
$$
\underset{\gamma\rightarrow \frac{1}{p}}\lim\frac{\sin\pi p\gamma}{\pi(1-p\gamma)}\Bigl(\frac{\pi}{\sin\pi\gamma}\Bigl)^p \frac{\delta^{1-p \gamma}}{p \gamma}=\Biggl(\frac{\pi}{\sin\frac{\pi}{p}}\Biggl)^p,
$$
we get
$$
\|S\|_{L^p\rightarrow L^p}\geq \frac{\pi}{\sin\frac{\pi}{p}}.
$$
To prove that $S$ is non compact,  assume  the contrary. Then $\lim_{\mathrm{mes}(D)\rightarrow 0}\|P_D S\|_{L^p\to L^p} =0$, where $P_Df:=\chi_Df$  ($\chi_{D}$ denotes  the characteristic function of the subset $D\subset [0,1]$) \cite[Theorem 3.1]{KZ}.

On the other hand, let $D_a:=[a, 1]$ and $x_a:=(1/\mathrm{mes}(D_a)^{1/p})\chi_{D_a}$.
Then $\|x_a\|^p_{L^p}=1$ and
$$
Sx_a=\frac{1}{\mathrm{mes}(D_a)^{1/p}}\int_{D_a}\frac{dt}{1-tz}=\frac{1}{(1-a)^{1/p}}\int_a^1\frac{dt}{1-tz}\geq
$$
$$
\geq\frac{1}{(1-a)^{1/p}}\frac{1}{1-az}\int_a^1 dt=(1-a)^{1-1/p}\frac{1}{1-az}.
$$
Therefore
$$
\|P_{D_a}Sx_a\|_{L^p}^p=\int_0^1 |\chi_{D_a}(z)|^p|Sx_a(z)|^p dz\geq (1-a)^{p-1}\int_a^1\frac{dz}{(1-az)^p}=
$$
$$
=\frac{1}{p-1}\frac{1}{a}\Biggl(1-\frac{1}{(1+a)^{p-1}}\Biggl).
$$
It follows that
$$
\underset{a\rightarrow 1}\limsup  \|P_{D_a}S\|_{L^p\to L^p}^p\geq \underset{a\rightarrow 1}\limsup \|P_{D_a}Sx_a\|_{L^p}^p\geq \frac{1}{p-1}\left(1-\frac{1}{2^{p-1}}\right)>0,
$$
which is a contradiction. This completes the proof on noncompactness of the operator $S$.

To show that the Markov-Stieltjes  transform  on $L^2 (0,1)$ is unitarily equivalent to   the Markov-Stieltjes  transform on  $H^2 (\mathbb{D})$, consider the restriction operator
$$
j : H^2(\mathbb{D}) \rightarrow L^2(0,1), f \mapsto f|(0,1).
$$
Let $S_L$ denotes the  Markov-Stieltjes  transform  on $L^2 (0,1)$, and $S_H$  the  Markov-Stieltjes  transform  on  $H^2(\mathbb{D})$. Note that $j S_H=S_L j$, i.e.,
 $j$ is an intertwining operator for $S_L$ and  $S_H$. The operator $j$ is bounded (see the formula (4)) and injective  and has a dense range.  By the Putnam-Douglas Theorem (see, e.g., \cite[Theorem IX.6.10(c)]{Con}),  $S_L$ is unitarily equivalent to  $S_H$. Application of Theorem 4.2 completes the proof.

\begin{corollary}
The Markov-Stieltjes transform  $S$ is a bounded operator from $H^p$ to $L^p(0,1)$ and
$$
\|S\|_{H^p \rightarrow L^p}\leq \pi^{1+1/p}\cot \frac{\pi}{2\max\{p,q\}}.
$$
\end{corollary}

 Proof.
By the formula (4), if $f\in H^p(\mathbb{D})$ then  the restriction $f|(0,1)$ belongs to
$L^p (0,1)$ and the norm of the restriction operator
$j_p\colon H^p \rightarrow L^p(0,1)\colon f\mapsto f| (0,1)$
does not exceed $\pi^{1/p}$. Let  $S_{H}$ denotes the Markov-Stieltjes
operator in  $H^p$. Then $S=j_pS_{H}: H^p  \rightarrow L^p (0,1)$ and
$$
\|S\|_{H^p \rightarrow L^p}\leq\|j_p\|_{H^p \to L^p}
\|S_{H}\|_{H^p\rightarrow H^p}\leq \pi^{1+1/p}\cot \frac{\pi}{2\max\{p,q\}}.
$$

\begin{corollary}
The map $f\mapsto L^pR_n(Sf)$ is continuous on  $L^p$.
\end{corollary}

 Proof. The proof is similar to the case of $H^p$.

\begin{remark}
The Markov-Stieltjes  transform   is an unbounded  operator on $L^1(0,1)$ and $L^\infty(0,1)$, but it is bounded as an operator from  $L^p(0,1)\ (1<p<\infty)$ to $L^1(0,1)$  \cite[p. 187]{Mir}. It is also a continuous map from $L^1(0,1)$  to $L^p(0,1)$ for $0<p<1$ because
 if $f,g\in L^p(0,1)$ then
$$
\|S(f-g)\|_{L^p(0,1)}^p=\int_0^1\Biggl|\int_0^1\frac{f(t)-g(t)}{1-tz}dt\Biggl|^p dz\leq\int_0^1\Biggl(\int_0^1\frac{|f(t)-g(t)|}{1-tz}dt\Biggl)^p dz
\leq
$$
$$
\leq\int_0^1\Biggl(\int_0^1\frac{|f(t)-g(t)|}{1-z}dt\Biggl)^p dz=\int_0^1\frac{dz}{(1-z)^p}\|f-g\|_{L^1(0,1)}^p=\frac{1}{1-p}\|f-g\|_{L^1(0,1)}^p.
$$
\end{remark}

Recall that the  Banach space
$\ell^p_A\ (1<p\leq \infty)$ consists of functions $f$ that are holomorphic on the unit disc,   $f(z)=\sum_{n=0}^\infty f_n z^n  (z\in \mathbb{D})$, and such that $\|f\|_{\ell^p_A}^p:=\sum_{n=0}^\infty |f_n|^p<\infty$. Obviously, the space $\ell^p_A$ can be identified with $\ell^p$. According to \cite[p. 53]{Bott}, a Hankel operator $H(a)$ on $\ell^p_A$
associated with a sequence $a=(a_n)$ is defined by
\[
H(a)f(z)=\sum_{j=0}^\infty b_j z^j\qquad  (z\in \mathbb{D}),
\]
where
\[
b_j=\sum_{k=0}^\infty a_{k+j+1}f_k\qquad \left(f(z)=\sum_{n=0}^\infty f_n z^n\right).
\]

 In the following theorem we consider $S$ as an operator on  $\ell^p_A$ (this means that  $S$ is defined by the formula (1), where $f\in \ell^p_A, z\in \mathbb{D}$).

\begin{theorem}
 The Markov-Stieltjes  transform $S$ is a bounded   Hankel operator on  $\ell^p_A\ (1<p<\infty)$ and has Hilbert matrix with respect to  the standard basis. Moreover,
\[
\| S\|_{\ell^p_A\rightarrow \ell^p_A}= \frac{\pi}{\sin \frac{\pi}{p}}.
\]
 \end{theorem}

 Proof.
 First note, that the Markov-Stieltjes  transform exists for $f\in \ell^p_A$, since  $f|(0,1)\in L^1(0,1)$. In fact, if $f(t)=\sum_{n=0}^\infty f_n t^n$ then
$$
\int_0^1|f(t)|dt\leq\int_0^1\sum_{n=0}^\infty |f_n|t^n dt=\sum_{n=0}^\infty\frac{|f_n|}{n+1},
$$
the series converges by the H$\ddot{\rm o}$lder inequality.

Arguing as in the proof of Theorem 4.2 we get that the Markov-Stieltjes  transform $S$ is a Hankel operator on  $\ell^p_A$ and has Hilbert matrix with respect to  the standard basis $\{e_n:n\in \mathbb{Z}_+\},\ e_n(z):=z^n$ of  $\ell^p_A$.

To compute the norm of $S$, note that the Hardy-Littlewood-Polya-Shur inequality \cite[Theorem 318]{HLP} applied to the function $K(x,y)=1/(x+y)$ implies, that
\[
\|Sf\|_{\ell^p_A}^p=\sum_{n=0}^\infty\left|\sum_{m=0}^\infty\frac{f_m}{n+m+1}\right|^p\leq k^p\|f\|_{\ell^p_A}^p,
\]
where the best possible constant $k$ is  \cite[p. 229]{HLP}
\[
k=\int_{0}^\infty\frac{dx}{x^{1/p}(1+x)}=\frac{\pi}{\sin\frac{\pi}{p}}
\]
(for the last equality see, e.g., \cite[Section 2.2.4, formula 25]{BMP}).

\begin{corollary}
Let $1<p\leq 2$.  The Markov-Stieltjes  transform $S$  is a bounded  operator from  $\ell^p_A$ to $H^q\ (1/p+1/q=1)$ and
  $$
  \|S\|_{\ell^p_A\rightarrow H^q}\leq\frac{\pi}{\sin\frac{\pi}{p}}.
  $$
\end{corollary}

 Proof.
It follows from the above theorem, because by  \cite[Theorem 6.1, p. 94]{Dur}  $\ell^p_A\subset H^q(\mathbb{D})$  and the norm of the natural embedding of $\ell^p_A$ into $H^q(\mathbb{D})$  does not exceed $1$.

\begin{corollary}
The Markov-Stieltjes transform  $S$ is a bounded operator from $\ell^p_A$ to $L^q(0,1)$ and
$$
\|S\|_{\ell^p_A \to L^q}\leq \frac{\pi^{1+1/q}}{\sin\frac{\pi}{q}}.
$$
\end{corollary}

 Proof. The proof follows from formula (4) and  Corollary 5.7.

 \end{document}